\documentclass[11pt]{amsart}
\usepackage[headings]{fullpage}
\usepackage{amssymb,epic,eepic,epsfig,amsbsy,amsmath,amscd,color}
\usepackage[all]{xy}
\usepackage{graphicx,subcaption}
\usepackage{texdraw}
\usepackage{url}
\usepackage{bbm}
\usepackage{mathrsfs}
\usepackage{standalone}
\usepackage{tikz}
\usetikzlibrary{cd,positioning,knots,patterns,hobby}
\usepackage{accents}
\usepackage[normalem]{ulem}
\usepackage{comment}

\theoremstyle{plain}
\newtheorem{theorem}{Theorem}[section]

\newtheorem{lemma}[theorem]{Lemma}
\newtheorem{corollary}[theorem]{Corollary}

\theoremstyle{definition}

\newtheorem{remark}[theorem]{Remark}

\newcommand\FIGc[3]{\begin{figure}[htpb]
    \includegraphics[height=#3]{#1.eps}
    \caption{#2}
    \label{fig:#1}
    \end{figure}}

\newcommand\incleps[2]{{\includegraphics[height=#1]{#2.eps}}}
\def\bdD{{\mathbf D}}
\def\cD{\mathcal D}

\DeclareMathOperator{\order}{ord}
\DeclareMathOperator{\MaxSpec}{MaxSpec}
\DeclareMathOperator{\Irr}{Irreps}
\DeclareMathOperator{\Hom}{Hom}

\newcommand{\marked}{\mathcal{P}}
\newcommand{\stateS}{\mathscr{S}}

\newcommand{\matchedS}{\mathscr{S}^{ma}}
\newcommand{\evenS}{\mathscr{S}^{ev}}

\def\BC{\mathbb C}
\def\BN{\mathbb N}
\def\BZ{\mathbb Z}

\def\BT{\mathbb T}

\def\cR{{\mathcal R}}

\def\la{\langle}
\def\ra{\rangle}

\def\al{\alpha}
\def\ve{\varepsilon}

\def\bD{{\bar \Delta }}

\newcommand\no[1]{}

\def\bP{\bar P}

\def\bT{\mathbb T}

\def\rk{\mathrm{rk}}
\def\cM{\mathcal M}

\def\cS{\mathscr S}
\def\ot{\otimes}

\def\cE{\mathcal E}

\def\bk{\mathbf k}

\def\cP{\mathcal P}

\def\fS{\mathfrak S}

\def\ev{{\mathrm{ev}}}
\def\cX{\mathcal X}

\def\embed{\hookrightarrow}

\def\cH{\mathcal H}

\def\bQ{\bar Q}
\def\ev{{\mathrm{ev}}}
\def\cX{\mathcal X}
\def\cX{\mathcal X}

\def\embed{\hookrightarrow}

\def\SM{(\Sigma,\cP)}
\def\pS{\partial \Sigma}

\def\pM{\partial M}
\def\cN{\mathcal N}
\def\MN{(M,\cN)}

\def\fr{\operatorname{fr}}

\def\cB{\mathcal B}

\def\ord{\mathrm{ord}}
\def\Sx{\cS_\xi}

\def\Sx{\cS_\xi}

\def\MNp{(M',\cN')}

\def\cSp{\cS^+}
\def\cD{\mathcal D}

\begin{document}

\title[Skein algebra]{Stated skein modules of marked 3-manifolds/surfaces, a survey}

\author[Thang  T. Q. L\^e]{Thang  T. Q. L\^e}
\address{School of Mathematics, 686 Cherry Street,
 Georgia Tech, Atlanta, GA 30332, USA}
\email{letu@math.gatech.edu}
\author[Tao Yu]{Tao Yu}
\address{School of Mathematics, 686 Cherry Street,
 Georgia Tech, Atlanta, GA 30332, USA}
\email{tyu70@gatech.edu}


\thanks{
2010 {\em Mathematics Classification:} Primary 57N10. Secondary 57M25.\\
}

\begin{abstract}
We give a survey of some old and new results about the stated skein modules/algebras of 3-manifolds/surfaces. For generic quantum parameter, we discuss the splitting homomorphism for the 3-manifold case, general structures of the stated skein algebras of marked surfaces (or bordered punctured surfaces) and their embeddings into quantum tori. For roots of 1 quantum parameter, we discuss the Frobenius homomorphism (for both marked 3-manifolds and marked surfaces), describe the center of the skein algebra of marked surfaces, the dimension of the skein algebra over the center, and the representation theory of the skein algebra. In particular, we show that the skein algebra of non-closed marked surface at any root of 1 is a maximal order. We give a full description of the Azumaya locus of the skein algebra of the puncture torus and give partial results for closed surfaces.
\end{abstract}
\maketitle

\tableofcontents

\section{Introduction}

We survey some old and new results about the stated skein modules of marked 3-manifolds and skein algebras of surfaces. 

Kauffman bracket skein modules for 3-manifolds and skein algebras for surfaces were introduced by Przytycki \cite{Przy} and Turaev \cite{Turaev,Turaev2} in around 1987. The skein modules/algebras have played an important role in low dimensional topology and quantum algebra as they have applications and connections to objects such as character varieties \cite{Bullock,PS,Turaev,BFK}, the Jones polynomial and its related topological quantum field theory (TQFT) \cite{Kauffman, BHMV,Tu:Book}, (quantum) Teichm\"uller spaces and (quantum) cluster algebras \cite{BW0,FGo, FST, Muller}, the AJ conjecture \cite{FGL,Le5}, and many more. 

More than 20 years after the Kauffman bracket skein algebra was introduced, Bonahon and Wong \cite{BW0} made an important contribution, proving that the skein algebra of a surface with at least one puncture can be embedded into a quantum torus by the quantum trace map, which is a quantization of the map expressing the trace of curves in the shear coordinates of Teichm\"uller space. The work opens possibilities to quantize Thurston's theory of hyperbolic surfaces to build hyperbolic TQFT and to better geometrically understand the volume conjecture \cite{Kashaev1}. One main problem is to understand the representation theory of the skein algebras at roots of 1.

Bonahon and Wong's proof of the existence of the quantum trace map suggests that the skein algebra of a surface can be split into smaller simple blocks. The first author \cite{Le:triangular} worked out and made precise this splitting phenomenon by introducing the stated skein algebra for punctured bordered surfaces, which are punctured surfaces having boundary. The main feature is the existence of the splitting homomorphism relating the stated skein algebra of a surface and that of its splitting along an ideal arc. This approach gives a new proof of the existence of the quantum trace map, and more \cite{CL}. The stated skein algebra theory fits well with the {\em integral} quantum group associated to $SL_2(\BC)$ and its integral dual $\mathcal O_{q^2}(\mathfrak{sl}_2)$, and many algebraic facts concerning the quantum groups have simple transparent interpretations by geometric formulas. For example, the stated skein algebra of the bigon $\cB(\cD_2)$, with its natural cobraided structure, is isomorphic to the cobraided Hopf algebra $\mathcal O_{q^2}(\mathfrak {sl}_2)$, and under the isomorphism the natural basis of $\cB(\cD_2)$ maps to Kashiwara's canonical basis of $\mathcal O_{q^2}(\mathfrak{sl}_2)$, see \cite{CL} and Section \ref{sec.surfaces}.

In this paper, we extend the definition of the stated skein modules to marked 3-manifolds (Section \ref{sec.2}), and then survey some old and new results on the stated skein modules/algebras. In Section \ref{sec.splitting} we sketch a proof of the splitting homomorphism for the 3-manifold case. Section \ref{sec.surfaces} gives general properties of skein algebras of surfaces: basis, Noetherian domain, Gelfand-Kirillov dimension. We also discuss relations between stated skein algebras and both quantum groups (for the ideal bigon and the ideal triangle) and the quantum moduli algebra of Alekseev-Grosse-Schomerus and Buffenoir-Roche (for once-marked surfaces).
Section \ref{sec.embedding} explains two embeddings of the stated skein algebra of a surface with at least one ideal point into quantum tori and relations between them. For closed surfaces we discuss the embedding of an associated graded algebra into a quantum torus. Section \ref{sec.Che} gives the Chebyshev-Frobenius homomorphism for marked 3-manifolds. In Section \ref{sec.root1} we describe the center of the stated skein algebra of a marked surface at a root $\xi$ of 1, calculate its PI degree, discuss the Azumaya locus in the general case, and give a precise description of the Azumaya locus for the punctured torus and partial results for closed surfaces. Proofs of some results are sketched, and details of many results will appear elsewhere.

\subsection{Acknowledgments}
The authors would like to thank F. Bonahon, F. Costantino, C. Frohman, J. Kania-Bartoszynska, T. Schedler, A. Sikora, and M. Yakimov for helpful discussions. The first author is supported in part by NSF grant DMS 1811114.

\section{Skein modules/algebras} 
\label{sec.2}

\subsection{Ground ring, notations}
Throughout the paper, the ground ring $\cR$ is a commutative Noetherian domain with a distinguished invertible element $q^{1/2}$. Let $\BN, \BZ, \BC$ denote respectively the sets of non-negative integers, integers, and complex numbers.

\subsection{Skein modules of marked 3-manifolds}
As in \cite{Le:qtrace}, a {\em marked 3-manifold} is a pair $(M,\cN)$, where $M$ is an oriented 3-manifold with (possibly empty) boundary $\partial M$, and $\cN \subset \partial M$ is a 1-dimensional oriented submanifold, called the {\em marking}, such that each connected component of $\cN$ is diffeomorphic to the closed interval $[0,1]$. 

An {\em $\cN$-tangle in $M$} is a compact 1-dimensional non-oriented submanifold $\al$ of $M$, equipped with a non-tangent vector field, called the {\em framing}, such that $\partial \al = \al \cap \cN$ and the framing at each boundary point of $\al$ is a positive tangent vector of $\cN$. Two $\cN$-tangles are {\em $\cN$-isotopic} if they are isotopic in the class of $\cN$-tangles. The empty set is considered as a $\cN$-tangle isotopic only to itself.

A {\em stated $\cN$-tangle} $\al$ is a $\cN$-tangle equipped with a  map $s : \partial \al \to \{\pm \}$, called the states of $\al$. 

\begin{figure}[h]
\begin{subfigure}[b]{0.45\linewidth}
\centering
\[
\raisebox{-0.43\height}{
\begin{tikzpicture}[scale=0.75]
\clip (0,0)circle(1);
\fill[gray!20!white] (0,0)circle(1);
\begin{knot}[clip width=5,background color=gray!20!white]
\strand (1,1)--(-1,-1);
\strand (-1,1)--(1,-1);
\end{knot}
\end{tikzpicture}
}
=
q\raisebox{-0.43\height}{
\begin{tikzpicture}[scale=0.75]
\clip (0,0)circle(1);
\fill[gray!20!white] (0,0)circle(1);
\draw (-1,-1)..controls (0,0)..(-1,1);
\draw (1,-1)..controls (0,0)..(1,1);
\end{tikzpicture}
}
+
q^{-1}\raisebox{-0.43\height}{
\begin{tikzpicture}[scale=0.75]
\clip (0,0)circle(1);
\fill[gray!20!white] (0,0)circle(1);
\draw (-1,-1)..controls (0,0)..(1,-1);
\draw (-1,1)..controls (0,0)..(1,1);
\end{tikzpicture}
}
\]
\subcaption{Skein relation}
\end{subfigure}
\hfill
\begin{subfigure}[b]{0.45\linewidth}
\centering
\[
\raisebox{-0.43\height}{
\begin{tikzpicture}[scale=0.75]
\clip (0,0)circle(1);
\fill[gray!20!white] (0,0)circle(1);
\draw (0,0)circle(.5);
\end{tikzpicture}
}
=
(-q^2-q^{-2})\raisebox{-0.43\height}{
\begin{tikzpicture}[scale=0.75]
\clip (0,0)circle(1);
\fill[gray!20!white] (0,0)circle(1);
\end{tikzpicture}
}
\]
\subcaption{Trivial knot relation}
\end{subfigure}
\begin{subfigure}[b]{0.45\linewidth}
\centering
\[
\raisebox{-0.57\height}[0.5\height][0.4\height]{
\begin{tikzpicture}
\fill[gray!20!white] (0,0)circle(.75);
\begin{scope}
\clip (0,0)circle(.75);
\begin{knot}[clip width=5,background color=gray!20!white,end tolerance=1pt]
\strand (-.15,-.9) to[in angle=0,
	curve through={(0,-.75) (.3,-0.4)}] (0,0);
\strand (.15,-.9) to[in angle=180,
	curve through={(0,-.75) (-.3,-0.4)}] (0,0);
\end{knot}
\end{scope}
\node at (.45,-0.6){$+$};
\node at (-.45,-0.6){$-$};
\fill (0,-.75) node[below]{$\mathcal{N}$} circle(2pt);
\end{tikzpicture}
}
=q^{-1/2}
\raisebox{-0.57\height}[0.5\height][0.4\height]{
\begin{tikzpicture}
\fill[gray!20!white] (0,0)circle(.75);
\fill (0,-.75) node[below]{$\mathcal{N}$} circle(2pt);
\end{tikzpicture}
}
\]
\subcaption{Returning arc relation 1}
\end{subfigure}
\hfill
\begin{subfigure}[b]{0.45\linewidth}
\centering
\[
\raisebox{-0.57\height}[0.5\height][0.4\height]{
\begin{tikzpicture}
\fill[gray!20!white] (0,0)circle(.75);
\begin{scope}
\clip (0,0)circle(.75);
\begin{knot}[clip width=5,background color=gray!20!white,end tolerance=1pt]
\strand (-.15,-.9) to[in angle=0,
	curve through={(0,-.75) (.3,-0.4)}] (0,0);
\strand (.15,-.9) to[in angle=180,
	curve through={(0,-.75) (-.3,-0.4)}] (0,0);
\end{knot}
\end{scope}
\node at (.45,-0.6){$+$};
\node at (-.45,-0.6){$+$};
\fill (0,-.75) node[below]{$\mathcal{N}$} circle(2pt);
\end{tikzpicture}
}
=0=
\raisebox{-0.57\height}[0.5\height][0.4\height]{
\begin{tikzpicture}
\fill[gray!20!white] (0,0)circle(.75);
\begin{scope}
\clip (0,0)circle(.75);
\begin{knot}[clip width=5,background color=gray!20!white,end tolerance=1pt]
\strand (-.15,-.9) to[in angle=0,
	curve through={(0,-.75) (.3,-0.4)}] (0,0);
\strand (.15,-.9) to[in angle=180,
	curve through={(0,-.75) (-.3,-0.4)}] (0,0);
\end{knot}
\end{scope}
\node at (.45,-0.6){$-$};
\node at (-.45,-0.6){$-$};
\fill (0,-.75) node[below]{$\mathcal{N}$} circle(2pt);
\end{tikzpicture}
}
\]
\subcaption{Returning arc relation 2}
\end{subfigure}
\begin{subfigure}[b]{\linewidth}
\centering
\[
\raisebox{-0.57\height}[0.5\height][0.4\height]{
\begin{tikzpicture}
\fill[gray!20!white] (0,0)circle(.75);
\begin{scope}
\clip (0,0)circle(.75);
\begin{knot}[clip width=5,background color=gray!20!white,end tolerance=1pt]
\strand (-.15,-.9)--(0,-.75)--(30:.75);
\strand (.15,-.9)--(0,-.75)--(150:.75);
\end{knot}
\end{scope}
\node at (0.4,-0.6){$-$};
\node at (-0.4,-0.6){$+$};
\fill (0,-.75) node[below]{$\mathcal{N}$} circle(2pt);
\end{tikzpicture}
}
=q^2
\raisebox{-0.57\height}[0.5\height][0.4\height]{
\begin{tikzpicture}
\fill[gray!20!white] (0,0)circle(.75);
\begin{scope}
\clip (0,0)circle(.75);
\begin{knot}[clip width=5,background color=gray!20!white,end tolerance=1pt]
\strand (-.15,-.9)--(0,-.75)--(30:.75);
\strand (.15,-.9)--(0,-.75)--(150:.75);
\end{knot}
\end{scope}
\node at (0.4,-0.6){$+$};
\node at (-0.4,-0.6){$-$};
\fill (0,-.75) node[below]{$\mathcal{N}$} circle(2pt);
\end{tikzpicture}
}
+q^{-1/2}
\raisebox{-0.57\height}[0.5\height][0.4\height]{
\begin{tikzpicture}
\fill[gray!20!white] (0,0)circle(.75);
\draw (30:.75)..controls (0,-.75)..(150:.75);
\fill (0,-.75) node[below]{$\mathcal{N}$} circle(2pt);
\end{tikzpicture}
}
\]
\subcaption{State exchange relation}
\end{subfigure}
\caption{Defining Relations in $\stateS(M,\mathcal{N})$}
\label{stated-relations}
\end{figure}
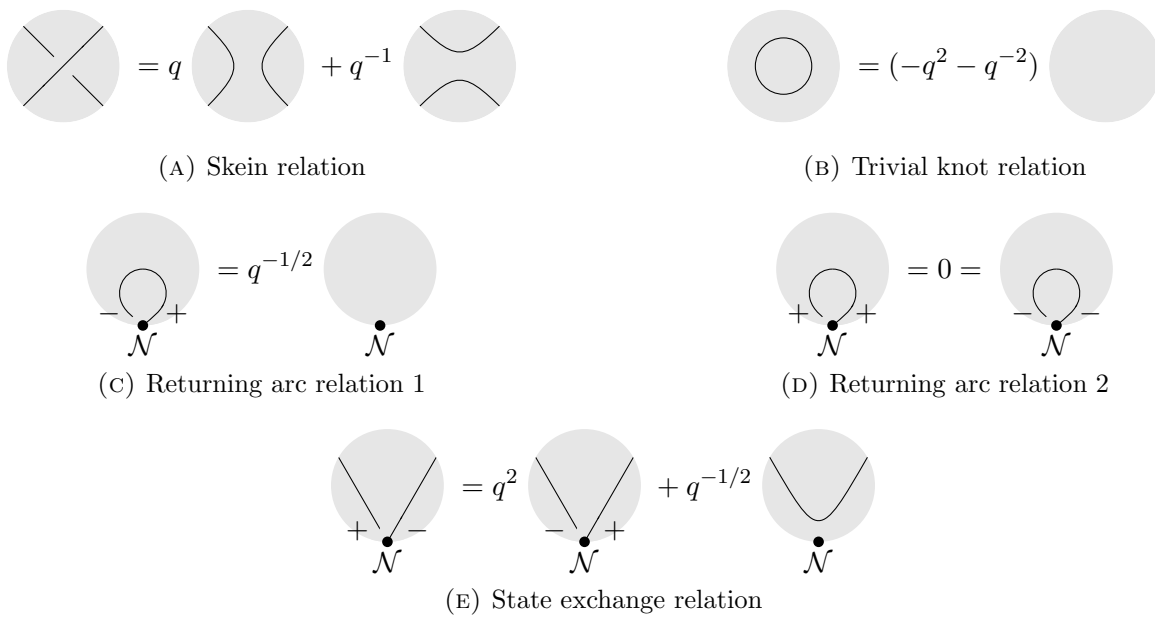

The {\em skein module} $\cS\MN$ is the $\cR$-module freely spanned by $\cN$-isotopy classes of $\cN$-tangles modulo the local relations (A-E) described in Figure \ref{stated-relations}. In each relation, the diagrams represent parts of $\cN$-tangles. The framings in the diagrams are perpendicular to the page and pointing to the reader. In (C-E), we assume that $\cN$ is perpendicular to the page, and its intersection with the page is the bullet labeled by $\cN$ there. There are two strands of the $\cN$-tangle coming to $\cN$, the lower one being depicted by the broken line. There are no other strands ending on the segment of $\cN$ between the two strands. The states of the endpoints are marked by $\pm$.

\begin{remark}
When $\cN=\emptyset$ we don't need the relations (C-E). In this case the skein module was introduced independently by J. Przytycki \cite{Przy} and V. Turaev \cite{Turaev,Turaev2}. Relations (A) and (B) were introduced by Kauffman \cite{Kauffman} in his study of the Jones polynomial. Relations (C) and (D), for surface case, appeared in \cite{BW0}. Relation (E) appeared in \cite{Le:triangular}.
\end{remark}

\subsection{Marked surfaces and punctured bordered surfaces}

A {\em marked surface} is a pair $\SM$, where $\Sigma$ is a compact oriented surface with possibly empty boundary $\pS$  and $\cP\subset \pS$ is a finite set, called the set of marked points. The  {\em associated marked 3-manifold} $\MN$ is defined by $M = \Sigma \times [-1,1]$ and $\cN=\cP \times [-1,1]$. Although $M$ has corners when $\pS\neq \emptyset$, we can smooth the corner and consider $M$ as a smooth 3-manifold. Define $\cS\SM= \cS \MN$ as an $\cR$-module. Given two stated $\cN$-tangles $\alpha, \al'$ define the product $\al \al'$ by stacking $\al$ above $\al'$. This gives $\cS\SM$ an $\cR$-algebra structure, which was first introduced by Turaev \cite{Turaev} for the case $\cP=\emptyset$ in connection with the quantization of the Atiyah-Bott-Weil-Petersson-Goldman symplectic structure of the character variety. The algebra $\cS\SM$ is closely related to the character variety (see \cite{Bullock,PS,BFK}) and the quantum Teichm\"uller space (see \cite{CF,Kashaev,BW0}).

\def\cR{{\mathcal R}}
\def\pfS{\partial \fS}
\def\tpfS{\widetilde{\partial \fS}}
 
The skein algebra $\cS\SM$ coincides with the stated skein algebra introduced by the first author in \cite{Le:triangular}, although at first they look different. Let us explain the connection here.



A boundary component of $\Sigma$ is {\em unmarked} if it does not have a marked point. Let $\fS$ be the result of removing all marked points and unmarked boundary components from $\Sigma$. Such a surface $\fS$ is called a {\em punctured bordered surface}. Each marked point of $\Sigma$ is called a {\em boundary ideal point} of $\fS$, and each unmarked boundary component of $\Sigma$ is called an {\em interior ideal point} of $\fS$. By adding back the ideal points one can recover $\SM$ from $\fS$. Each connected component $c$ of the boundary $\pfS$ is an open interval, called a {\em boundary edge} of $\fS$, and $c\times (-1,1)$ is called a {\em boundary wall} of $M':= \fS \times (-1,1)$. The boundary $\tpfS$ of $M'$ is the disjoint union of all the boundary walls.

By a {\em $\tpfS$-tangle} $\al$ in $M'= \fS\times (-1,1)$ we mean a framed 1-dimensional compact non-oriented submanifold properly embedded in $M'$ with vertical framing at each endpoint and distinct heights for endpoints in each boundary wall. Two $\tpfS$-tangles are {\em $\tpfS$-isotopic} if they are isotopic in the class of $\tpfS$-tangles. Note that the endpoints of $\al$ in one boundary wall are linearly ordered by heights since they have distinct heights, and $\tpfS$-isotopy does not change the height order.

Every $\tpfS$-tangle can be represented by a {\em $\pfS$-tangle diagram}, which by definition\footnote{In this paper a $\pfS$-tangle diagram is a positively ordered $\pfS$-tangle diagram of \cite{Le:triangular}} is a tangle diagram on $\fS$ where
\begin{itemize}
\item the endpoints of the tangle diagram are distinct points in $\pfS$, and
\item on each boundary edge $c$ the height order is given by the positive direction of $c$ (inherited from the orientation of $\fS$).
\end{itemize}

Every $\cN$-tangle in $M$ can be turned into a $\tpfS$-tangle in $M'$ by slightly moving the endpoints above a marked point on a boundary component $b$ along the negative direction of $b$ (opposite of the arrows in Figure~\ref{boundary-relations}). This gives a bijection of $\cN$-isotopy classes of $\cN$-tangles and $\tpfS$-isotopy classes of $\tpfS$-tangles. Hence the skein algebra $\cS\SM$ is canonically isomorphic to the skein algebra $\cS(\fS)$ defined as the $\cR$-module freely generated by $\tpfS$-isotopy classes of stated $\tpfS$-tangles modulo the original relations (A) and (B), and the new relations (C-E) of Figure \ref{boundary-relations}, which are the translations of the relations (C-E) in Figure \ref{boundary-relations}. This is the definition introduced in \cite{CL}.

\begin{figure}[h]
\begin{subfigure}[b]{0.45\linewidth}
\centering
\[
\raisebox{-0.4\height}{
\begin{tikzpicture}[scale=1.2]
\fill[gray!20!white] (0,0)rectangle(1,1);
\draw[arrows={-Stealth[scale=1.5]}] (1,0)--(1,1);
\draw (1,0.3)node[right]{$-$}
	..controls (0.3,0.3) and (0.3,0.6).. (1,0.6)node[right]{$+$};
\end{tikzpicture}
}
=q^{-1/2}
\raisebox{-0.4\height}{
\begin{tikzpicture}[scale=1.2]
\fill[gray!20!white] (0,0)rectangle(1,1);
\draw[arrows={-Stealth[scale=1.5]}] (1,0)--(1,1);
\end{tikzpicture}
}
\]
\phantomsubcaption
\phantomsubcaption
\subcaption{Returning arc relation 1}
\end{subfigure}
\hfill
\begin{subfigure}[b]{0.45\linewidth}
\centering
\[
\raisebox{-0.4\height}{
\begin{tikzpicture}[scale=1.2]
\fill[gray!20!white] (0,0)rectangle(1,1);
\draw[arrows={-Stealth[scale=1.5]}] (1,0)--(1,1);
\draw (1,0.3)node[right]{$+$}
	..controls (0.3,0.3) and (0.3,0.6).. (1,0.6)node[right]{$+$};
\end{tikzpicture}
}
=0=
\raisebox{-0.4\height}{
\begin{tikzpicture}[scale=1.2]
\fill[gray!20!white] (0,0)rectangle(1,1);
\draw[arrows={-Stealth[scale=1.5]}] (1,0)--(1,1);
\draw (1,0.3)node[right]{$-$}
	..controls (0.3,0.3) and (0.3,0.6).. (1,0.6)node[right]{$-$};
\end{tikzpicture}
}
\]
\subcaption{Returning arc relation 2}
\end{subfigure}
\begin{subfigure}[b]{\linewidth}
\centering
\[
\raisebox{-0.4\height}{
\begin{tikzpicture}[scale=1.2]
\fill[gray!20!white] (0,0)rectangle(1,1);
\draw[arrows={-Stealth[scale=1.5]}] (1,0)--(1,1);
\draw (0,0.3)--(1,0.3)node[right]{$+$};
\draw (0,0.6)--(1,0.6)node[right]{$-$};
\end{tikzpicture}
}
=q^2
\raisebox{-0.4\height}{
\begin{tikzpicture}[scale=1.2]
\fill[gray!20!white] (0,0)rectangle(1,1);
\draw[arrows={-Stealth[scale=1.5]}] (1,0)--(1,1);
\draw (0,0.3)--(1,0.3)node[right]{$-$};
\draw (0,0.6)--(1,0.6)node[right]{$+$};
\end{tikzpicture}
}
+q^{-1/2}
\raisebox{-0.4\height}{
\begin{tikzpicture}[scale=1.2]
\fill[gray!20!white] (0,0)rectangle(1,1);
\draw[arrows={-Stealth[scale=1.5]}] (1,0)--(1,1);
\draw (0,0.3)..controls (0.7,0.3) and (0.7,0.6)..(0,0.6);
\end{tikzpicture}
}
\]
\subcaption{State exchange relation}
\end{subfigure}
\caption{Boundary relations}
\label{boundary-relations}
\end{figure}

\subsection{Positive state submodule}

The submodule $\cSp\MN$ of $\cS\MN$ spanned by $\cN$-tangles with positive states was introduced in \cite{Le:qtrace}. The corresponding subalgebra $\cSp\SM\subset\stateS(\Sigma,\marked)$ in the  marked surface was first defined by Muller \cite{Muller} in connection with quantum cluster theory. The Muller algebra $\cSp\SM$ is a quantization of the decorated Teichm\"uller space of Penner \cite{Penner}.

\subsection{Category of marked 3-manifolds}

The skein module can be considered as a functor from the category of marked $3$-manifolds and (isotopy classes of) embeddings to the category of $\cR$-modules. 

An {\em embedding} of a marked 3-manifold $\MN$ into $(M',\cN')$ is an orientation preserving proper embedding $f: M \embed M'$ that restricts to an orientation preserving embedding on $\cN$. An embedding $f$ induces an $\cR$-module homomorphism $f_\ast: \cS\MNp \to \cS\MNp$ by $f_\ast[T]=[f(T)]$ for any stated $\cN$-tangle $T$.
Clearly $f_\ast$ depends only on the isotopy class of $f$. A {\em morphism} from $\MN$ to $(M',\cN')$ is an isotopy class of embeddings from $\MN$ to $(M',\cN')$.


If $\MN = (M_1, \cN_1) \sqcup (M_2, \cN_2)$, then there is a natural isomorphism
\begin{equation}
\cS (M_1, \cN_1) \ot_R \cS(M_2, \cN_2)\cong\cS\MN
\end{equation}
sending $[T_1]\otimes[T_2]$ to $[T_1\cup T_2]$ for any $\mathcal{N}_1$-tangle $T_1$ and $\mathcal{N}_2$-tangle $T_2$.

In the case of marked surfaces, there is a similar picture. The skein algebra is a monoidal functor from the category of marked surfaces and embeddings to the category of $\cR$-algebras. An {\em embedding} of a marked surface $\SM$ into $(\Sigma', \marked')$ is an orientation preserving proper embedding $f: \Sigma \embed \Sigma'$ such that  $f(\cP) \subset \cP'$. The embedding $f$ induces an $\cR$-algebra homomorphism $f_\ast: \cS\SM \to \cS(\Sigma', \cP')$. 

\subsection{Higher rank group $SL_n$}
The theory of skein modules/algebras for Lie group $SL_n$ has been developed in joint work of the first author and A. Sikora \cite{LS}, where analogs of a splitting homomorphism (Theorem \ref{thm.splitting}) and a quantum group isomorphism of the bigon algebra (Theorem \ref{thm.iso1}) are obtained. 

\def\tal{\tilde \al}

\section{Splitting homomorphism} \label{sec.splitting}

A basic property of the skein module of marked 3-manifolds is the existence of the splitting homomorphism which relates the skein module of a marked 3-manifold to a new, simpler, marked 3-manifold obtained by splitting the original one along an embedded disk. The splitting homomorphism has many applications and is one of the most important technical tools. It was in a search for such a splitting homomorphism that prompted the first author to introduce the stated skein algebras of marked surfaces. 

\subsection{Splitting homomorphism for 3-manifolds}
Suppose $\MN$ is a marked 3-manifold and $a_1,a_2$ are two connected components of the marking set $\cN$. We do not assume $M$ is connected.

For each $i=1,2$, choose a closed disk $D_i\subset \pM$ whose interior contains $a_i$ such that $D_1$ and $D_2$ are disjoint and disjoint from any other connected component of $\cN$. Each disk $D_i$ inherits an orientation from $M$. Choose an orientation reversing diffeomorphism $h: D_1 \to D_2$ such that $h$ restricts to an orientation preserving diffeomorphism $a_1\to a_2$. Let $M'= M/(D_1\equiv_h D_2)$ be the 3-manifold obtained from $M$ by identifying $D_1$ with $D_2$ via $h$, and let $\cN' = \cN \setminus (a_1 \cup a_2) \subset \partial M'$. Then $(M', \cN')$ is a marked 3-manifold (after smoothing the corner).

There is a natural projection $p: M \to M'$. Let $a=p(a_1)=p(a_2)$ and  $D=p(D_1)=p(D_2)$. Then $D$ is a disk properly embedded in $M'$, disjoint from $\cN'$, and containing the arc $a$ in its interior. We will say $\MN$ is the result of splitting $(M', \cN')$ along $(D,a)$.
 
An $\cN'$-tangle $\al$ in $M'$ is said to be {\em $(D,a)$-transversal} if $p^{-1}(\al)$ is an $\cN$-tangle, i.e.
\begin{itemize}
\item $\al$ is transversal to $D$, and $\al \cap D = \al \cap a$, and
\item the framing at every point of $\al \cap a$ is a positive tangent vector of $a$.
\end{itemize}
Then $\tal := p^{-1}(\al)$ is an $\cN$-tangle in $M$. If in addition $\al$ is stated, then $\tal$ is stated at every endpoint except for the endpoints on $a_1 \cup a_2$. Given a map 
$$s: \al \cap a \to \{\pm\},$$
let $(\tal, s)$ be the stated $\cN$-tangle whose state at a point $x \in \tal \cap (a_1 \cup a_2)$ is given by $s(p(x))$.

\begin{theorem} \label{thm.splitting}
Assume $\MN$ is the result of splitting $(M', \cN')$ along $(D,a)$, with the above notations. There is a unique $\cR$-module homomorphism
$$ \Theta: \cS(M',\cN') \to \cS\MN$$
such that if $\al$ is an $\cN'$-tangle $\al$ in $M'$ which is $(D,a)$-transversal, then
$$ \Theta(\al) = \sum_{s: \al \cap a \to \{\pm\}} (\tal,s).$$
\end{theorem}

\begin{remark}
If the splitting is applied to multiple disks, it can be done in different orders. It is clear from definition that the splitting homomorphism is independent of the order.
\end{remark}

\begin{proof}
We sketch a proof here. For details, see \cite{BL}. The majority of the work is done in \cite{Le:triangular}.  

Let $T(D)$ denote the $\cR$-module freely generated by $(D,a)$-transversal $\mathcal{N}'$-tangles (not isotopy classes). The formula above defines a map
\[\tilde{\Theta}:T(D)\to\stateS(M,\mathcal{N}).\]

We need to show that $\tilde{\Theta}$ is invariant under isotopy and all moves given by the defining relations (A-E) of Figure \ref{stated-relations}. We can assume that the supports of the isotopy and the moves are small. If the support is disjoint from $D$ then clearly $\tilde{\Theta}$ is invariant. On the other hand, if the support is small and intersects $D$, the invariance of $\tilde{\Theta}$ is verified in the proof of Theorem~3.1 of \cite{Le:triangular}.
\end{proof}

The spitting homomorphism gives the skein theory of marked 3-manifolds some flavor of a topological quantum field theory. We will explore this direction in the upcoming work \cite{CL2}.

\subsection{Splitting homomorphism for surfaces}
Suppose $p_1$ and $p_2$ are two marked points of a marked surface $\SM$ which might be disconnected. For $i=1,2$ choose a closed interval $b_i\subset \pS$ which contains $p_i$ in its interior. By identifying $b_1$ with $b_2$ via an orientation reversing diffeomorphism which maps $p_1$ to $p_2$, from $\Sigma$ we get a new surface $\Sigma'$. Let $\cP'= \cP\setminus \{p_1, p_2\}$, considered as a subset of $\pS'$. The splitting homomorphism gives an $\cR$-linear homomorphism
\begin{equation}
\Theta: \cS(\Sigma', \cP') \to \cS\SM. \label{eq.iso1}
\end{equation}
For surfaces, we get a stronger result.

\begin{theorem} \cite{CL} 
The splitting homomorphism in \eqref{eq.iso1} is an injective $\cR$-algebra homomorphism.
\end{theorem}

Originally, the splitting homomorphism is defined for punctured bordered surfaces. In this case, splitting a punctured bordered surface  $\mathfrak{S}'$ along an {\em ideal arc} $e$ gives a new punctured bordered surface $\fS$ with new boundary edges $e_1, e_2$ (among all boundary edges) such that after gluing $e_1$ with $e_2$ one recovers $\fS'$. The splitting homomorphism gives an algebra embedding
$\Theta: \cS(\fS') \to \cS(\fS)$.

\section{Surfaces} \label{sec.surfaces}
We discuss general properties of the skein algebra of a marked surface and its relations to known algebras, including the quantum groups associated to $SL_2(\BC)$ and their canonical bases.

Throughout $\SM$ is a marked surface, with corresponding punctured bordered surface $\fS$. Let
\begin{equation}
r\SM= r(\fS):=
\begin{cases}
0, &\text{if $\fS$ is the sphere with no or one ideal point}, \\
1, &\text{if $\fS$ is the sphere  with two ideal points},\\
2, &\text{if $\fS$ is the closed torus},\\
3 m - 3 \chi(\fS), &\text{otherwise}.
\end{cases} 
\label{eq.rfS}
\end{equation}
Here $m=|\cP|$ is the number of marked points or the number of boundary edges of $\fS$, and $\chi(\Sigma)=\chi(\fS)$ is the Euler characteristic.

\subsection{Basis, domain, Gelfand-Kirillov dimension}

Unlike the 3-manifold case, the skein module $\cS(\fS)$ of a punctured bordered surface is always a free $\cR$-module, with a basis described below.

A $\pfS$-tangle diagram $\al$ is {\em simple} if it has no crossing, no trivial loop, and no trivial arc. Here a {\em trivial loop} is a simple closed curve bounding a disk in $\fS$, and a {\em trivial arc} is one which can be homotoped relative to its endpoints into a boundary edge. 
 


\def\cR{{\mathcal R}}
We order the set $\{ \pm \}$ so that $+$ is greater than $-$. The state $s: \partial \al \to \{\pm\}$ of a $\pfS$-tangle diagram $\al$ is {\em increasing} if when traversing any boundary edge along its positive direction, the state values are never decreasing, i.e. one never encounters a $-$ immediately after a $+$.

Let $B(\fS)$ be the set of all isotopy classes of increasingly stated, simple $\pfS$-tangle diagrams, including the empty set, which by convention is considered as a simple $\pfS$-tangle diagram. 

\def\inv{\mathrm{nd}}
\begin{theorem}\label{thm.basis}
(a) \cite{Le:triangular} As an $\cR$-module, $\cS(\fS)$ is free with basis $B(\fS)$.

(b) \cite{LY} As an $\cR$-algebra, $\cS(\fS)$ is orderly finitely generated by arcs and loops. This means, there are one-component $\pfS$-tangle diagrams $x_1, \dots, x_n \in B(\fS)$ such that the set $\{ x_1^ {k_1} \dots x_n^{k_n} \mid k_i \in \BN\}$ spans $\cS(\fS)$ over $\cR$.

(c) \cite{LY} $\cS(\fS)$ is a Noetherian domain. 

(d) \cite{LY} The Gelfand-Kirillov dimension of $\cS(\fS)$ is $r(\fS)$.
\end{theorem}

When $\fS$ has no boundary, parts (a), (b), (c) were known, see respectively \cite{Przy,AF,PS2}. The finite generation (without order, and for no boundary $\fS$) was first proved in \cite{Bullock}. When $\fS$ has at least one ideal point and no boundary, the fact that $\cS(\fS)$ is a domain follows from \cite{BW0} where it is proved that the quantum trace map embeds $\cS(\fS)$ into a quantum torus, which is a domain. Part (d) follows from Theorems \ref{thm.embedding}(b) and \ref{thm.embedding2} below.
 
\def\letterx{  \raisebox{-9pt}{\incleps{.9 cm}{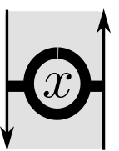}} } 
\def\lettery{  \raisebox{-9pt}{\incleps{.9 cm}{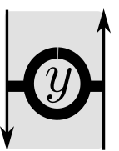}} }
\def\letterxy{  \raisebox{-9pt}{\incleps{.9 cm}{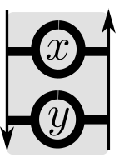}} }

\subsection{Bigon, co-braided quantum group, and canonical basis} 
The bigon $\bdD_2$ is the disk with two marked points on the boundary; its corresponding punctured bordered surface, denoted by $\cD_2$, is the disk with two boundary ideal points, see Figure~\ref{fig:bigon}.

\FIGc{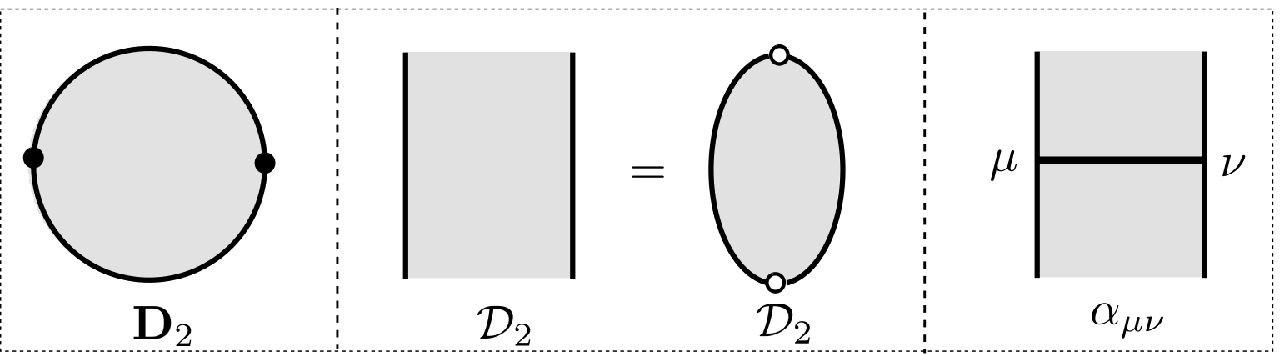}{Bigons: $\mathbf D_2$, $\cD_2$, and arcs $\alpha_{\mu\nu} $ in $\cD_2$ with $\mu,\nu \in \{\pm \}$.}{2.5cm}

The topology of the bigon $\cD_2$ is very special, which allows us to define many interesting operations on $\cS(\cD_2)$. For example, by splitting the bigon $\cD_2$ along an ideal arc connecting the two ideal points, we get a disjoint union of two bigons. The corresponding splitting algebra homomorphism
$$ \Delta: \cS(\cD_2) \to \cS(\cD_2) \ot_\cR \cS(\cD_2)$$
defines a coproduct on $\cS(\cD_2)$. Similarly both a counit and an antipode may be geometrically defined, making $\cS(\cD_2)$ a Hopf algebra, see \cite{CL}. Moreover, the Hopf algebra $\cS(\cD_2)$ is {\em cobraided}, i.e. it has a co-$R$-matrix $ \rho: \cS(\cD_2) \ot_\cR \cS(\cD_2) \to \cR$ which turns the category of $\cS(\cD_2)$-comodules into a braided category. For an overview of the theory of cobraided Hopf algebra see \cite{Majid}. The co-$R$-matrix is given by the following nice geometric formula, where $x$ and $y$ are $\partial \cD_2$-tangle diagrams and $\ve$ is the counit.
$$\rho \left(\letterx  \ot \lettery \right)  = \epsilon \left(  \letterxy \right).$$

It turns out that $\cS(\cD_2)$ is isomorphic to the well-known cobraided Hopf algebra $\mathcal{O}_{q^2}(\mathfrak{sl}_2)$, which is the Hopf dual of the quantized enveloping algebra $U_{q^2}(\mathfrak{sl}_2)$ (see for example \cite{Majid}). As an $\cR$-algebra $\mathcal{O}_{q^2}(\mathfrak{sl}_2)$ is generated by $a,b,c,d$ modulo the relations
\begin{align*}
&ca=q^2ac,\,db=q^2bd,\,ba=q^2ab,\,dc=q^2cd,\,bc=cb,\\
&ad-q^{-2}bc=da-q^2bc=1.
\end{align*}

\begin{theorem}\label{thm.iso1}\cite{Le:triangular,CL}
The four elements $\al_{\pm,\pm} \in \cS(\cD_2)$  in Figure \ref{fig:bigon} generate the $\cR$-algebra $\cS(\cD_2)$. The map given by
\[\alpha_{+,+} \mapsto a ,\  \alpha_{+,-} \mapsto b ,\ \alpha_{-,+} \mapsto c ,\ \alpha_{-,-} \mapsto d\] 
is an isomorphism from $\cS(\cD_2)$ to $\mathcal{O}_{q^2}(\mathfrak{sl}_2)$ as cobraided Hopf algebras. The basis $B(\cD_2)$, up to powers of $q$, maps to Kashiwara's canonical basis of $\mathcal{O}_{q^2}(\mathfrak{sl}_2)$.
\end{theorem} 

Kashiwara's canonical basis of $\mathcal{O}_{q^2}(\mathfrak{sl}_2)$ was one of the first constructed canonical bases of quantum groups \cite{Kashiwara}. The canonical basis has many nice properties, including the positivity. We will explore positivity of the basis $B(\fS)$ for other surfaces such as polygons, in a future work.

\subsection{Comodule algebras over $\mathcal{O}_{q^2}(\mathfrak{sl}_2)$} Let $c$ be a boundary edge of a punctured border surface $\fS$. By gluing $c$ to the left edge of $\cD_2$, we get a back $\fS$. The splitting homomorphism gives a map
$$ \Delta_c : \cS(\fS) \to \cS(\fS) \ot_\cR \cS(\cD_2).$$
It turns out that $\Delta_c$ gives $\cS(\fS)$ a right $\cS(\cD_2)$-comodule structure. Moreover the comodule structure is compatible with algebra structure of $\cS(\fS)$ and makes $\cS(\fS)$ a right $\cS(\cD_2)$-comodule algebra, as defined in \cite{Kassel,Majid}. 

Any comodule over $\mathcal{O}_{q^2}(\mathfrak{sl}_2)$ is a module over $U_{q^2}(\mathfrak{sl}_2)$. The structure of integrable $U_{q^2}(\mathfrak{sl}_2)$-modules for generic $q$ is well-known. In \cite{CL} it was proved that $\cS(\fS)$ is an integrable $U_{q^2}(\mathfrak{sl}_2)$-module with an explicit description of the highest weight vectors. The $\cS(\cD_2)$-comodule structure allows us to characterize the image of the splitting homomorphism \eqref{eq.iso1} in terms of the 0-th Hochschild homology, a result obtained independently by \cite{CL} and \cite{KQ}.

More generally given a marked 3-manifold $\MN$ and a component $c$ of the marking $\cN$, the splitting homomorphism also gives $\cS\MN$ a comodule structure over $\cS(\cD_2)$, and its image is closely related to the  0-th Hochschild homology.

\def\cD{\mathcal D}

\subsection{Attaching an ideal triangle}
Let $\cD_3$ be the ideal triangle with boundary edges $c_1, c_2, c_3$ as in Figure \ref{fig:triangle}.

\FIGc{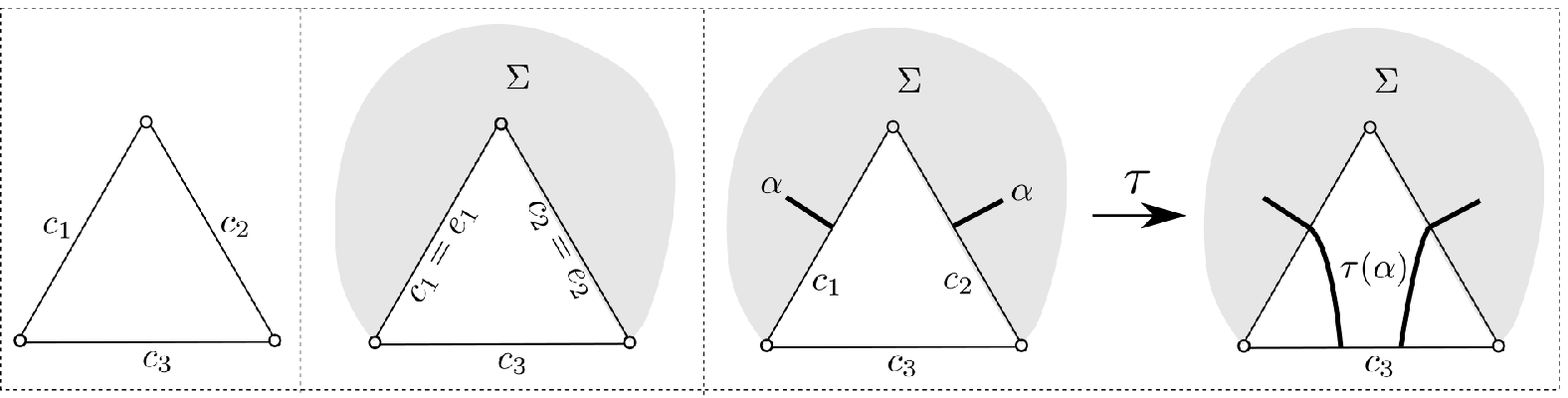}{Ideal triangle, attaching an ideal triangle to $\Sigma$, and the map $\tau$}{3cm}

Suppose $e_1, e_2$ are two boundary edges of a bordered puncture surface $\fS$. Let $\fS'$ be the result of attaching $\cD_3$ to $\fS$ by identifying $c_1$ with $e_1$ and $c_2$ with $e_2$. For a $\pfS$-tangle diagram $\al\in B(\fS)$ let $\tau(\al)$ be the $\pfS'$-tangle diagram obtained by extending the arcs ending in $c_1 \cup c_2$ to arcs ending in $c_3$ as in Figure \ref{fig:triangle}. Since $B(\fS)$ is an $\cR$-basis of $\cS(\fS)$, the map $\tau$ extends to an $\cR$-linear homomorphism $\tau: \cS(\fS) \to \cS(\fS')$, which is not an algebra homomorphism in general.

\begin{theorem}\cite{CL}\label{thm.glue2}
The map $\tau: \cS(\fS) \to \cS(\fS')$ is bijective.
\end{theorem}

The algebra structure of $\cS(\fS')$ can be described explicitly using that of $\cS(\fS)$ via the co-braiding structure of $\cS(\cD_2)$.  Let us only mention two special cases, both are taken from \cite{CL}.

First suppose $\fS =\fS_1 \sqcup \fS_2$ with $e_1\subset \fS_1$ and $e_2\subset \fS_2$. Then each $\cS(\fS_i)$, with boundary edge $e_i$, is a $\mathcal{O}_{q^2}(\mathfrak{sl}_2)$-comodule-algebra. Given any two $\mathcal{O}_{q^2}(\mathfrak{sl}_2)$-comodule-algebras, their braided tensor product is defined, which is also a $\mathcal{O}_{q^2}(\mathfrak{sl}_2)$-comodule-algebra, see \cite{Majid}. The explicit description of the algebra structure shows that $\cS(\fS')$ is the braided tensor product of $\cS(\fS_1)$ and $\cS(\fS_2)$.

In the second example $\fS=\cD_2$, the bigon, which has two boundary edges. Attaching an ideal triangle gives $\cM_1$, which is the once-punctured monogon, see Figure \ref{fig:monogon}. The above consideration gives a proof that $\cS(\cM_1)$ is isomorphic to the braided version (or transmutation) of $\mathcal{O}_{q^2}(\mathfrak{sl}_2)$, defined in \cite[Examples 4.3.4 and 10.3.3]{Majid}.

\FIGc{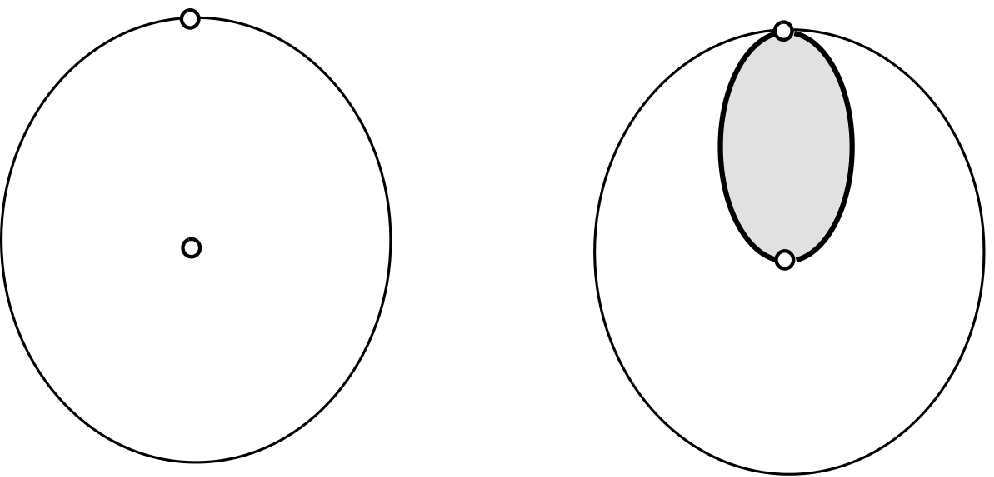}{Once-punctured monogon $\cM_1$ (left) as the result of attaching an ideal triangle  to a bigon (shaded).}{2.5cm}

More generally, by attaching ideal triangles we can convert the $n$-gon $\cD_n$ to an $(n-1)$-punctured monogon. The above procedure describes a relation between the two skein algebras.

\subsection{Quantum moduli algebra, factorization homology}
The results of the previous subsections allow to build $\cS(\fS)$ from a few simple algebras. In the special case when there is exactly one marked point this helps to identify $\cS(\fS)$ with known algebras. Assume the ground ring $\cR$ is a field.

Suppose $\Sigma_{g,n}$ is the compact surface of genus $g$ having $n$ boundary components. For each braided Hopf algebra $H$,  Alekseev, Grosse, and Schomerus \cite{AGS} and Buffenoir and  P. Roche \cite{BR} define the quantum moduli algebra $A(\Sigma_{g,n};H)$ which is a quantization (of Fock-Rosly's Poisson structure \cite{FR}) of the moduli space of flat $G$-connections of $\Sigma_{g,n}$ when $H$ is the quantum group associated with the simply-connected Lie group $G$. The algebra $A(\Sigma_{g,n};H)$ later appeared in the work of Ben-Zvi, Brochier, and Jordan \cite{BBJ}, where they showed that the factorization homology of $\Sigma_{g,n}$ with value in the category of $H$-modules is the category of $A(\Sigma_{g,n};H)$-modules in an appropriate sense.

\begin{theorem}\cite{LY} \label{thm.iso4}
For $H=U_{q^2}(\mathfrak{sl}_2)$ the quantum moduli algebra $A(\Sigma_{g,n};H)$ is isomorphic to the skein algebra $\cS(\Sigma_{g,n}')$, where $\Sigma_{g,n}'$ is the result of removing an open disk from $\Sigma_{g,n}$ and marking a point on the newly created boundary component. 
\end{theorem}
We learned that Theorem \ref{thm.iso4} was also independently proved by M. Faitg \cite{Faitg}.



\def\pbfS{\partial \bar{\fS }}
\section{Embedding into quantum tori}  \label{sec.embedding}

Quantum tori are a class of simple algebras with nice properties. In this section we discuss embeddings of skein algebras of marked surfaces into quantum tori. 

\subsection{Quantum tori} 
Informally, a quantum torus is an algebra of Laurent polynomials in several variables which $q$-commute, i.e. $ab= q^k ba$ for some integer $k$. By definition, the {\em quantum torus} associated to an anti-symmetric $r\times r$ integral matrix $Q$ is
\[\mathbb{T}(Q):=\cR\langle x_1^{\pm1},\dots,x_r^{\pm1}\rangle/\langle x_ix_j=q^{Q_{ij}}x_jx_i\rangle.\]

A quantum torus is a Noetherian domain \cite{GW}. In particular, it has a ring of fractions, which is a division algebra. The Gelfand-Kirillov dimension of $\bT(Q)$ is $r$.

For $\bk=(k_1,\dots, k_r)\in \BZ^r$ let
$$ x^\bk = q^{-\frac{1}{2} \sum_{i<j} Q_{ij} k_i k_j} x_1 ^{k_1} x_2^{k_2} \dots x_r ^{k_r}.$$
Then $\{ x^\bk \mid \bk \in \BZ^r\}$ is a free $\cR$-basis of $\bT(Q)$, and
\begin{equation}
x^\bk x ^{\bk'} = q ^{\frac 12 \la \bk, \bk'\ra_Q}  x^{\bk + \bk '},
\end{equation}
where
$$ \la \bk, \bk'\ra_Q = \sum_{1\le i, j \le r} Q_{ij} k_i k'_j.$$

If $\Lambda\subset \BZ^r$ is a submonoid, then the $\cR$-submodule
$ A(Q;\Lambda)\subset \bT(Q)$ spanned by $\{ x^\bk, \bk \in \Lambda \}$ is an $\cR$-subalgebra of $\bT(Q)$, called a {\em monomial subalgebra}. When $\Lambda= \BN^r$, the corresponding subalgebra is denoted by $\bT_+(Q)$.

For any positive integer $N$ it is easy to check that there is an algebra embedding, called $N$-th {\em Frobenius homomorphism},
\begin{equation}
\Phi_N: \bT(N^2 Q) \to \bT(Q), \text{ given by } \Phi_N(x_i)= x_i^{N}. \label{eq.Fro}
\end{equation}

More generally, suppose $Q'$ is another anti-symmetric $r'\times r'$ integral matrix such that $HQ' H^T= Q$, where $H$ is an $r\times r'$ matrix and $H^T$ is its transpose. Then the $\cR$-linear map  $\bT(Q)\to \bT(Q')$ given on the basis by $x^\bk \to x^{\bk H}$, where $\bk H$ is the product of the horizontal vector $\bk$ and the matrix $H$, is an algebra homomorphism, called a {\em multiplicatively linear homomorphism}.

\subsection{Non-closed surfaces} \label{sec.Q}
For a punctured bordered surface $\fS$ the number $r(\fS)$ is defined by~\eqref{eq.rfS}. Recall that if $\fS$ has an ideal triangulation, then $r(\fS)= 3m - 3\chi$, where $\chi$ is the Euler characteristic and $m$ is the number of boundary components.


\begin{theorem} \cite{LY} \label{thm.embedding}
Suppose $\fS$ is a punctured bordered surface which has at least one ideal point (interior or on the boundary).

(a) There exists an $\cR$-algebra embedding $\varphi: \cS(\fS) \embed \bT(\bQ)$, where $\bQ$ is an anti-symmetric integral matrix of size $r(\fS) \times r(\fS)$. Moreover, there is an $\BN$-filtration of $\cS(\fS)$ compatible with the algebra structure such that the associated graded algebra is isomorphic to a monomial subalgebra $A(\bQ;\Lambda)$, where $\Lambda\subset \BZ^{r(\fS)}$ is a submonoid of  rank $r(\fS)$. 

(b) Suppose in addition $\pfS\neq \emptyset$. There is an anti-symmetric integral matrix $\bP$ of size $r(\fS) \times r(\fS)$ and an $\cR$-algebra embedding $\phi:\cS(\fS) \embed \bT(\bP)$ such that $\bT_+(\bP)\subset\phi( \cS(\fS)) \subset \bT(\bP)$. 
\end{theorem}

We give the definitions of $\bQ$ and $\bP$ below, which depend on ideal (quasi-)triangulations of $\fS$. We also give a sketch of a proof of part (b).

The relation between $\varphi$ and $\phi$ is as follows. There is a bigger quantum torus $\bT(\bP')$ containing $\bT(\bP)$ as a monomial subalgebra and a multiplicatively linear algebra homomorphism $\Psi:\bT(\bQ)\to \bT(\bP')$ which maps $\varphi$ to $\phi$. Loosely speaking, when $\pfS=\emptyset$, the quantum torus $\bT(\bP)$ can be considered as a quantization of the decorated Teichm\"uller space \cite{Penner} using Penner's lambda length coordinates, while $\bT(\bQ)$ can be considered as a quantization of the enhanced Teichm\"uller space \cite{BW0} (or holed Teichm\"uller space of \cite{FG}) using shear coordinates. The map $\Psi$  is a quantization of  the map changing Penner's coordinates to shear coordinates. The quantum  Teichm\"uller spaces was first defined in \cite{CF,Kashaev}.  

One has $\rk(\bQ)=\rk(\bP)= r\SM-b_\ev$, where $b_\ev$ is the number of boundary components of $\Sigma$ having an even number (including 0) of marked points. 

\begin{remark} (i) When $\fS$ has no boundary, Theorem \ref{thm.embedding}(a) was first proved by Bonahon and Wong \cite{BW0}, and was an important development in the theory. Later other proofs were given in \cite{Le:qtrace,Le:triangular}. The map $\varphi$ is called the quantum trace map.

(ii) When there are no interior ideal points, the restriction of $\phi$ to the subalgebra $\cS^+(\fS)$ was first constructed by Muller \cite{Muller}. The extension to the case when $\fS$ has interior ideal points, again for the subalgebra $\cS^+(\fS)$, was done in \cite{LP}. When $\fS$ has no boundary edges the connection map $\Psi$ was first constructed by the first author \cite{Le:qtrace}.
\end{remark}

\subsection{Description of the matrix $\bQ$ of Theorem \ref{thm.embedding}} \label{sec.Q1}
  
\def\cornerab{  \raisebox{-9pt}{\incleps{.9 cm}{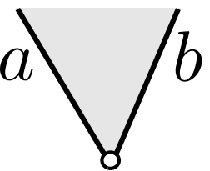}} } 
\def\cornerba{  \raisebox{-9pt}{\incleps{.9 cm}{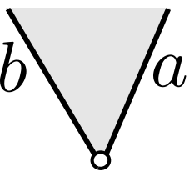}} }
\def\Dd{\Delta_\partial}

Let us describe the matrix $\bQ$ (of Theorem \ref{thm.embedding}) when $\fS$ is not the sphere with one or two ideal points. This matrix depends on an ideal triangulation of $\fS$, which we recall first.

An {\em ideal arc} of $\fS$ is a proper embedding of $(0,1)$ into $\fS$. Under a proper embedding, the ideal points of $(0,1)$ map to (not necessarily distinct) ideal points of $\fS$.
An {\em ideal triangulation} $\Delta$ of $\fS$ is a maximal collection of disjoint, pairwise non-isotopic essential ideal arcs. Here an arc is essential if it does not bounds a disk in $\fS$. Then $\Delta$ must contain the set $\Dd$ of all boundary edges (up to isotopy).
By splitting $\fS$ along the non-boundary arcs in $\Delta$ we get a collection of ideal triangles.

Let $\Dd'=\{ e' \mid e\in \Dd\}$ be another copy of $\Dd$, and $\bD= \Delta\sqcup \Dd'$. Then $|\bD| = r(\fS)$.
Define the anti-symmetric matrix $\bQ: \bD \times \bD \to \BZ$ by
\begin{align}
\bQ(a',a)& = 2 = - \bQ(a,a'),& &\text{if } a\in \Dd,\notag \\
\bQ(a,b)& = \#\left(\cornerba\right)- \#\left(\cornerab\right),& &\text{if } a,b\in \Delta,\label{eq.Q}\\
\bQ(a,b)& = 0,                &    &\text{otherwise}\notag.
\end{align}
Here each shaded part is a corner of an ideal triangle. Thus, the right hand side of \eqref{eq.Q} is the number of corners where $b$ is clockwise to $a$ minus the number of corners where $a$ is clockwise to $b$. The restriction of $\bQ$ to $\Delta$, denoted by $Q$, is a well-known matrix in the theory of Teichm\"uller space, describing the Poisson structure of the enhanced (or holed) Teichm\"uller space in shear coordinates, see e.g. \cite{FG}. 

\def\bcE{\bar{\cE}}
\subsection{Matrix $\bP$ and sketch of a proof of Theorem \ref{thm.embedding}(b)}\label{sec.P}

Let us sketch a proof of Theorem~\ref{thm.embedding}(b). Assume $\SM$ is not the monogon or the bigon. First we  describe the matrix $\bP$, which is defined only when $\SM$ has at least one marked point. 

\begin{figure}[h]
\centering
\begin{subfigure}[b]{0.4\linewidth}
\centering
\begin{tikzpicture}
\fill[gray!20!white] (1.2,0)--(1.2,1.3)--(-1.2,1.3)--(-1.2,0);
\begin{scope}
\clip (-1.2,0)rectangle(1.2,1.3);
\begin{knot}[clip width=10,background color=gray!20!white,end tolerance=1pt]
\strand[very thick] (-0.1,-0.1)--(0,0)--(1,1.3);
\strand[very thick] (-1,1.3)--(0,0)--(0.1,-0.1);
\end{knot}
\end{scope}
\draw (-0.4,0.2)node{$+$} (0.4,0.2)node{$+$};
\draw (-1.2,0.5)node[left]{
};
\draw (-1.2,0)--(1.2,0);
\fill (0,0) circle(2pt);
\end{tikzpicture}
\subcaption{$e$ with identical endpoints}\label{fig-norm-plus}
\end{subfigure}
\begin{subfigure}[b]{0.4\linewidth}
\centering
\begin{tikzpicture}
\fill[gray!20!white] (1.5,0)--(1.5,1.3)--(-1.5,1.3)--(-1.5,0);
\draw (-1.5,0)--(1.5,0);
\fill (-1,0) circle(2pt) (1,0) circle(2pt);
\draw[very thick] (-1,0) edge node[above]{$\hat{e}$}
	node[above,pos=0.1]{$+$} node[above,pos=0.9]{$-$}(1,0);
\end{tikzpicture}
\subcaption{States of $\hat{e}$}\label{fig-e-hat}
\end{subfigure}
\caption{Height convention, states of boundary arcs $\hat e$}
\end{figure}
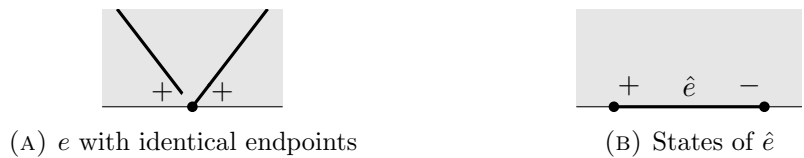

A \emph{$\marked$-arc} is a map $e:[0,1]\to\Sigma$ which embeds $(0,1)$  into $\Sigma\setminus\marked$, and $e(0), e(1) \in \cP$. A $\marked$-arc $e$ defines an element, also denoted by $e$, of the skein algebra $\cS\SM$ by assigning $+$ states to both endpoints and using the vertical framing. This is well defined if the arc has two distinct endpoints. If the endpoints coincide, we fix the relative heights of the endpoints so that the right strand is higher as in Figure~\ref{fig-norm-plus}.

If a $\marked$-arc $e$ is contained in the boundary, it defines another element $\hat{e}$ in the skein algebra $\cS\SM$ by assigning states as in Figure~\ref{fig-e-hat}. The element $\hat e$ is called a bad arc in \cite{CL}.

\def\Ed{\cE_\partial}
\def\hEd{\widehat{\Ed}}
\def\bE{\bar {\cE}}
\def\Pab{  \raisebox{-8pt}{\incleps{.9 cm}{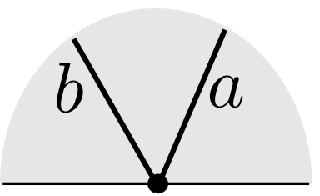}} }
\def\Pba{  \raisebox{-8pt}{\incleps{.9 cm}{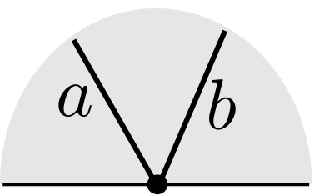}} }  

A \emph{quasitriangulation} $\cE$ is a maximal collection of disjoint, pairwise non-isotopic essential $\marked$-arcs. 
The non-boundary arcs in a quasitriangulation cut the surface into triangles and once-punctured monogons. Here a once-puncture monogon is an annulus with exactly one marked point. The collection $\cE$ must contain the set $\Ed$ of all boundary $\cP$-arcs. Denote $\hEd= \{ \hat e \mid e\in \Ed\}$. 

Let $\bcE= \cE \sqcup \hEd\sqcup \cH$, where $\cH$ is the set of all non-marked components of $\pS$. One can check that $|\cE|=r(\fS)=3 |\cP|- 3 \chi(\Sigma)$. Note that $\bE$ is a subset of the $\cR$-basis $B(\fS)$ of the skein algebra $\cS(\fS)= \cS\SM$, see Theorem \ref{thm.basis}(a). It is easy to verify that elements of $\bE$ are $q$-commuting, i.e. there is an anti-symmetric function $\bP: \bcE \times \bcE \to \BZ$ such that $xy = q^{\bP(x,y)} yx$ for $x,y\in \bE$. Using the defining relations of $\cS\SM$, $\bP$ can be described combinatorially. The result is (undefined entries are inferred by anti-symmetry)
\begin{align}
\bP(a,b) &=  0,  &      &\text{if } \ a \in \cH, b\in \bcE, \notag\\
\bP(a,b) &= \#\left( \Pab  \right)   - \# \left( \Pba  \right),  &    &\text{if} \ a, b\in \cE,  \label{eq.P1}\\
\bP(a,\hat b) &=  -\#\left( \Pab  \right)   - \#\left( \Pba  \right),   &    &\text{if} \ a\in \cE, b\in \Ed, \label{eq.P2}\\
\bP(\hat a, \hat b) &=  - \bP(a,b), &    &\text{if} \ a,b \in \Ed. \notag
\end{align}

Here the shaded parts are parts of $\Sigma$, with boundary the horizontal lines.
The $\cP$-arcs $a$ and $b$ are isotoped so that their interiors are in the interior of $\Sigma$ (even when one or both of them are boundary $\cP$-arcs). They are not necessarily distinct or edges of the same triangle, i.e. between them there might be other $\cP$-arcs.  

A monomial in the elements of $\bE$ is an element of the basis set $B(\fS)$ up to a scalar. Hence one sees that $\bT_+(\bP)\embed \cS\SM=\cS(\fS)$. One can show further that each monomial is a non-zero divisor, and for each $x\in \cS\SM$ there is a monomial $a$ such that $ax \in \bT_+$. From here one can prove that there is an embedding of $\cS\SM$ into $\bT(\bP)$.

The restriction $P$ of $\bP$ on the set $\cE\cup \cH$ was defined by Muller.

\def\rd{{\mathrm{rd}}}
\def\cSd{\cS^\rd}
\def\fM{\mathfrak{M}}
\subsection{Reduced skein algebra and quantum cluster algebra} Let $I$ be the ideal generated by all the elements $\hat e$ of Subsection \ref{sec.P}. The quotient $\cSd\SM= \cS\SM/I$ is called the reduced skein algebra in \cite{CL}, where it is proved that there is an embedding
\begin{equation}
\varphi^\rd :\cSd(\fS) \embed \BT(Q). \label{eq.em5}
\end{equation}

This embedding is essentially a quotient of $\varphi:\cS(\fS) \embed \BT(\bQ)$ given by Theorem \ref{thm.embedding}(a). In turn, the more general embedding $\varphi$ can be constructed through the embedding $\varphi^\rd$ of a bigger surface.

In \cite{Le:triangular}, it is proved that there is an embedding $\cS^+\SM\embed \BT(Q)$. As both $\cSd(\fS)$ and $\cS^+\SM$ embed into the same quantum torus $\BT(Q)$, one might want to compare them. It turns out that $\cSd\SM$ is an Ore localization of $\cS^+\SM$, which in turn is equal to the {\em quantum cluster algebra of $\SM$} when each boundary component has at least two marked points.
 
Let $\fM$ be the multiplicative set generated by all boundary $\cP$-arcs (with $+$ states) in $\cS^+\SM$. Muller showed in \cite{Muller} that $\fM$ is a two sided Ore set and the Ore localization $\cS^+\SM \fM^{-1}$ is equal to the quantum cluster algebra (as defined in \cite{BZ}) which quantizes the cluster algebra associated with the marked surface $\SM$.

\begin{theorem} \cite{LY}\label{thm.reduced-cluster}
The quantum cluster algebra $\cS^+\SM \fM^{-1}$ is naturally isomorphic to the reduced skein algebra $\cSd\SM$. 
\end{theorem}

Thus the reduced skein algebra gives a geometric model for the quantum cluster algebra.

When $\cP\neq \emptyset$,  $\cSd\SM$ and $\cS^+\SM$ embed into $\BT(P)$. Both matrices $P$ and $Q$ have  size $2m - 3 \chi(\Sigma)$ and nullity equal to the number of boundary components of $\Sigma$.

\subsection{Closed surfaces}
Suppose $\fS=\Sigma_g$ is a closed surface of genus $g$ (no boundary, no ideal points). These surfaces are exactly the ones not covered by Theorem \ref{thm.embedding}. 
 
When $g=1$ the algebra $\cS(\Sigma_1)$ embeds into a quantum torus \cite{FG}. However for $g\ge 2$ there is no known embedding of $\cS(\Sigma_g)$ into a quantum torus. The best we have so far is the following.
 
\begin{theorem} \label{thm.embedding2} \cite{KL}
Suppose $g\ge 2$. There is a $(6g-6)\times (6g-6)$ anti-symmetric integral matrix $Q$ with determinant $2^{6g-6}$ and an $\BN$-filtration of the algebra $\cS(\Sigma_g)$ such that the associated algebra embeds into the quantum torus $\bT(Q)$. Moreover, there is an $\BN \times \BZ$-filtration of the algebra $\cS(\Sigma_g)$, where $\BN\times \BZ$ is given the lexicographic order, and a $(6g-6)\times (6g-6)$ anti-symmetric integral matrix $Q'$ such that the associated algebra is isomorphic to $\bT_+(Q')$.
\end{theorem}
 
Here is an explicit description of $Q$. Let $\Sigma_{g', n'}$ be a surface of genus $g'$ with $n'$ ideal points, where $n' \ge 1$ and $2g'+n'-1=g$.
(Then the double of $\Sigma_{g', n'}$ along its ideal points gives $\Sigma_g$.) An ideal triangulation of $\Sigma_{g', n'}$ will have $3g-3$ edges, and gives rise to an anti-symmetric matrix $Q( \Sigma_{g', n'})$ of size $3g-3$ as described in Subsection \ref{sec.Q1}. Let $Q$ be anti-symmetric of size $6g-6$ given by
$$ Q = \begin{pmatrix}
Q( \Sigma_{g', n'}) & 2I \\
- 2I & 0 \end{pmatrix},$$
where each block is a square matrix of size $3g-3$, and $I$ is the identity matrix. The matrix $Q$ depends on a triangulation of $\Sigma_{g', n'}$. The doubles of the edges of the ideal triangulation of $\Sigma_{g', n'}$ give a pair of pants decomposition of $\Sigma_g$. The $\BN$-filtration of $\cS(\Sigma_g)$ of Theorem \ref{thm.embedding2} is based on this pair of pants decomposition, and is closely related to the one in \cite{PS2}.

A. Sikora also announced that he found a degeneration of $\cS(\Sigma_g)$ related to a quantum torus.

\section{The Frobenius homomorphism}
\label{sec.Che}

\subsection{Roots of 1}

When $\cR=\BC$ and $q=\xi$, a non-zero complex number, we denote $\cS\MN$ by $\Sx\MN$. Technically we have to fix also a square root of $\xi$, but the choice of such a square root does not play an important role in what follows.

For a root $\xi$ of 1 let $\order(\xi)$ be the least positive integer $n$ such that $\xi^n=1$.

\subsection{The 3-dimensional manifold case}
Let $\MN$ be a marked 3-manifold. A 1-component $\cN$-tangle $\al$ is diffeomorphic to either the circle $S^1$ or the closed interval $[0,1]$; we call $\al$ an {\em $\cN$-knot} in the first case and {\em $\cN$-arc} in the second case. For a 1-component $\cN$-tangle $\al$ and an integer $k\ge 0$, write $\al^{(k)} \in \cS_\xi(M)$ for the {\em $k$th framed power of $\al$} obtained by stacking $k$ copies of $\al$ in a small neighborhood of $\al$ along the direction of the framing of $\al$. Given a polynomial $P(z) = \sum c_i z^i \in \BZ[z]$, the {\em threading} of $\al$ by $P$ is given by $P^{\fr}(\al) = \sum c_i \al^{(i)} \in \cS_\xi(M)$.

The Chebyshev polynomials of type one $T_n(z) \in \BZ[z]$ are defined recursively as
\begin{equation}\label{eq.Che}
T_0(z)=2,\quad T_1(z)=z,\quad T_n(z) = zT_{n-1}(z)-T_{n-2}(z),\quad\forall n \geq 2.
\end{equation}

\begin{theorem}[See \cite{BL}]  \label{thm.1}
Suppose $(M,\cN)$ is a marked 3-manifold and $\xi$ is a complex root of unity. Let $N=\order(\xi^4)$ and $\ve =\xi^{N^2}$. 
 
There exists a unique $\BC$-linear map $\Phi_\xi: \cS_\ve(M,\cN) \to \cS_\xi(M,\cN)$ such that for any $\cN$-tangle $T = a_1 \cup \cdots \cup a_k \cup \al_1 \cup \cdots \cup \al_l$ where the $a_i$ are $\cN$-arcs and the $\al_i$ are $\cN$-knots,
\begin{align*}
\Phi_\xi(T) & = a_1^{(N)} \cup \cdots \cup a_k^{(N)} \cup (T_N)^{\fr}(\al_1) \cup \cdots \cup (T_N)^{\fr}(\al_l) \quad \text{in }\Sx\MN\\
&:= \sum_{0\le j_1, \dots, j_l\le N} c_{j_1} \dots c_{j_l}  a_1^{(N)} \cup \dots \cup  a_k^{(N)} \cup \, \al_1^{(j_1)} \cup \cdots \cup \al_l^{(j_l)} \quad \text{in }\Sx\MN,
\end{align*}
where $T_N(z) = \sum_ {j=0}^N c_j z^j$.

The map $\Phi_\xi$ is functorial and commutes with all the splitting homomorphisms.
\end{theorem}

We call $\Phi_\xi$ the {\em Chebyshev-Frobenius homomorphism}. 

Note that if $\order(\xi^4)=N$ then $\xi^{2N} =1$ or $\xi^{2N}=-1$. More precisely, $\xi^{2N} =(-1)^{N'+1}$, where $N'= \ord(\xi^2)$.
An addendum of Theorem~\ref{thm.1} is the following fact, which says the image of $\Phi_\xi$ is ``almost transparent" in the following sense.

\begin{theorem} \label{r.transparent}
Suppose $\MN$ is a marked 3-manifold, $\xi$ is a root of unity, $N'=\order(\xi^2)$. Then the image of $\Phi_\xi$ is almost transparent in $\cS_\xi\MN$ in the sense that
$$ \raisebox{-14pt}{\incleps{1.2 cm}{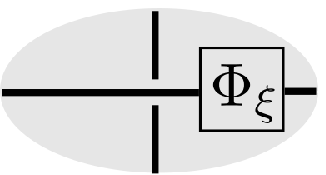}}  =  (-1)^{N'+1}\, \raisebox{-14pt}{\incleps{1.2 cm}{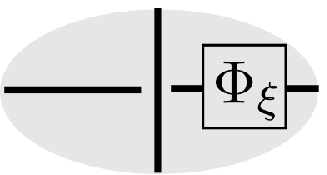}}.$$
Here the box enclosing $\Phi_\xi$ means one applies $\Phi_\xi$ to the component of the $\cN$-tangle containing the horizontal line, while the vertical line belongs to another component of the $\cN$-tangle.
\end{theorem}

\subsection{The case of surfaces}

For surfaces, the above results can be strengthened.

\begin{theorem}\cite{BL} \label{thm-frob-surface}
Suppose $\fS$ is a punctured bordered surface
and $\xi$ is a root of unity. Let $N=\order(\xi^4)$ and $\ve =\xi^{N^2}$. 

(a) There exists an algebra embedding  $\Phi_\xi: \cS_\ve(\fS) \to \cS_\xi(\fS)$ characterized by $\Phi_\xi(a)=a^{(N)}$ for all stated arcs $a$ on $\fS$, and $\Phi_\xi(\alpha)=T_N(\alpha)$ for all knots $\al$ on $\fS$.

(b) Assume $\fS$ has at least one ideal point. Let $Q$ and $\varphi$ be as in Theorem \ref{thm.embedding}(a). The following  diagram is commutative
\begin{equation}
\begin{tikzcd}
\cS_\ve \SM  \arrow[hookrightarrow, "\varphi"]{r} \arrow[d,"\Phi_\xi"]{d}& \bT(N^2Q)  \arrow[d, "\Phi_N"] \\
\cS_\xi \SM \arrow[hookrightarrow,"\varphi"]{r} & \bT(Q)
\end{tikzcd}
\label{eq.dia}
\end{equation}
\end{theorem}
Here $\bT(Q)$ and $\bT(N^2Q)$ are quantum tori defined over $\BC$ with $q=\zeta$, and $\Phi_N$ is the Frobenius homomorphism between quantum tori defined by \eqref{eq.Fro}. The upper right corner of the diagram should be the quantum torus of $Q$ defined with $q=\ve$, but it is equal to $\bT(N^2Q)$ since $\ve= \xi^{N^2}$. Diagram \eqref{eq.dia} shows that under the embeddings $\varphi$ into quantum tori, the Chebyshev-Frobenius becomes the usual Frobenius algebra homomorphism between quantum tori.

It should be noted that in general, if $a$ is a stated arc on $\fS$, we have
\begin{equation}
\Phi_\xi(a) = a^{(N)} \neq a^N, \label{eq.diff}
\end{equation}
even when $N$ is odd.

From Theorem \ref{r.transparent} one can deduce the following.
\begin{corollary} \label{r.center1}
Suppose $\fS$ is a punctured bordered surface
and $\xi$ is a root of unity, with $N=\order(\xi^4)$. If $\al$ is a simple closed curve on $\fS$, then $T_{2N}(\al)$ is in the center of $\cS(\fS)$. If $a\subset \fS$ is a stated $\pfS$-arc, then $a^{(4N)}$ is in the center of $\cS(\fS)$.
\end{corollary}


\begin{remark} (i) When $\pfS=\emptyset$, Theorem \ref{thm-frob-surface} was first proved by \cite{BW1}, using the quantum trace map. Another proof using skein modules was later given in \cite{Le:Frobenius}. 

(ii) For the submodule $\cS^+$ of positive states, Theorems \ref{thm.1}-\ref{thm-frob-surface} were proved in \cite{LP}.

(iii) In \cite{KQ}, part (a) of Theorem \ref{thm-frob-surface} is proved  for the case when $\ord(\xi)$ is odd. But even in this case, the definition of $\Phi_\xi$ in \cite{KQ} is different from ours: if $a$ is an arc, then in \cite{KQ} one has $\Phi_\xi(a)= a^N$, whereas we have $\Phi_\xi(a)= a^{(N)}$, which is different, see \eqref{eq.diff}.

(iv) Suppose $\xi$ is an arbitrary non-zero complex number, $N$ is an arbitrary positive integer, and $\ve= \xi^{N^2}$. In Diagram \eqref{eq.dia} all the maps are defined except for $\Phi_\xi$. One might ask for what values of $\xi$ and $N$ can the Frobenius map $\Phi_N$ restrict to a map $\Phi_\xi$ so that Diagram \eqref{eq.dia} commutes, and in addition, that $\Phi_\xi$ is functorial. The answer is that this happens exactly when $\xi$ is a root of 1 and $N=\ord(\xi^4)$. The proof for the similar statement for $\cS^+(\fS)$ was given in \cite{LP}, and this result can be used to prove the statement in our case.
\end{remark}

\subsection{Relation with Lusztig's quantum Frobenius homomorphism} For the bigon $\cS_{\xi}(\cD_2) = \mathcal O_{\xi^2}(\mathfrak{sl}_2)$,
 which is the Hopf dual of  (Lusztig's version of) the quantum group $U_{\xi^2}(\mathfrak{sl}_2)$.
Our Chebyshev-Frobenius  homomorphism $\Phi_\xi:  \cS_{\ve}(\cD_2) \to \cS_{\xi}(\cD_2)$, a Hopf algebra homomorphism, gives rise to a dual map, also a Hopf algebra homomorphism,
\begin{equation}
\Phi_\xi^\ast: U_{\xi^2}(\mathfrak{sl}_2) \to U_{\ve^2}(\mathfrak{sl}_2).
\end{equation}
Note that $\ve^2=\pm1$. When $\ve^2=1$, a simple check (see \cite{BL}) on generators shows that $\Phi^\ast_\xi$ is exactly the quantum Frobenius homomorphism of Lusztig \cite{Lusztig:1,Lusztig}. (Lusztig also mentioned the existence of a quantum Frobenius homomorphism for the case $\ve^2=-1$ although he did not write it down explicitly.) Thus our Chebyshev-Frobenius homomorphism can be interpreted as an extension of Lusztig's quantum Frobenius homomorphism from the bigon to all surfaces and marked 3-manifolds.

\section{Skein algebra at root of 1} \label{sec.root1}

The center of an algebra is important, for example, in questions about the representations of the algebra, see a review in Subsection \ref{sec.Azu}. We discuss the structure of skein algebra at roots of 1, describing its center, its PI dimension, and its representations. For the punctured torus, we give a full description of the Azumaya locus. We also extend a result of \cite{Jordan}, which states that the Azumaya locus of the skein algebra of a closed surface in the case at odd roots of 1 contains the smooth locus, to all cases of roots of 1.

When $q=\xi$ is a root of unity and there is no marked point, i.e.  $\cP=\emptyset$, the center of $\Sx(\Sigma, \emptyset)$ is determined in \cite{FKL} and is instrumental in proving the main result there, the unicity conjecture of Bonahon and Wong. Representations of skein algebras of non-marked surfaces at roots of 1 were initiated in \cite{BW0}--\cite{BW4}.

Throughout $\SM$ is a marked surface which is not $\Sigma_{0,0}, \Sigma_{0,1}, \Sigma_{0,2}$ and $\Sigma_{1,0}$ (with no marked points). As in Subsection \ref{sec.P}, we consider each $\cP$-arc $e$ as an element of $\cS\SM$ by assigning state $+$ to both endpoints. When $e$ is a boundary $\cP$-arc we also defined $\hat e$ as in Figure \ref{fig-norm-plus}.
 

\def\tZ{{\tilde Z}}
\def\tA{{\tilde A}}
\def\cC{\mathcal C}
\subsection{Representations of algebras and Azumaya locus, a review}
\label{sec.Azu}

Let $A$ be an associative $\BC$-algebra such that its center $Z$ is a domain with field of fractions $\tZ$. Let $\tA:= \tZ \ot_Z A$.  The {\em dimension of $A$ over $Z$}, denoted by $ \dim_ZA$, is the dimension of the vector space $\tA$ over the field $\tZ$.  

We make the following assumptions:
\begin{itemize}
\item[(*1)] $A$ is finitely generated as a $\BC$-algebra and  $A$ is a domain.
\item[(*2)] $A$ is finitely generated as a module over its center $Z$. 
\end{itemize}

By the Artin-Tate lemma, $Z$ is finitely generated as a $\BC$-algebra. Hence its maximal spectrum $\MaxSpec(Z)$ is an irreducible affine variety. It is known that $\tA$ is a division algebra with center $\tZ$. It follows that $\dim_\tZ\tA$ is a perfect square, $\dim_ZA= d^2$. The number $d$ is known as the {\em PI degree} of $A$.

Let $\Irr(A)$ denote the set of all equivalence classes of irreducible finite dimensional representations of $A$ over $\BC$. Since $Z$ is commutative, $\Irr(Z)$ is the set of all 1-dimensional representations,
$$ \Irr(Z) = \Hom_{\BC-\text{alg}} (Z, \BC) \overset \kappa\equiv  \MaxSpec(Z),$$
where $\Hom_{\BC-\text{alg}} (Z, \BC)$ is the set of all $\BC$-algebra homomorphisms from $A$ to $\BC$, and $\kappa(f) = \ker f$. By Schur's lemma, the restriction of an irreducible representation of $A$ to its center $Z$ gives a 1-dimensional representation of $Z$. Hence we have the central character map
\begin{equation}
\cC: \Irr(A) \to \Irr(Z) \equiv \MaxSpec(Z).
\end{equation}
The following summarizes some of the main facts about representations of $A$. 
\begin{theorem} (See e. g. \cite{BG,BY,DP}) \label{thm.Azu}
(a) Every irreducible representation of $A$ has dimension $\le d$, the PI dimension of $A$.

(b) The image under $\cC$ of all irreducible representations of dimension $d$ is a Zariski open dense subset $U$ of $\MaxSpec(Z)$. Moreover, every point of $U$ has exactly one preimage. Every point of $\MaxSpec(Z)$ has a non-zero finite number of preimages.
\end{theorem}

The set $U$ is known as the Azumaya locus of $A$. 

\subsection{Center at generic $q$} 
Let $\cH$ be the set of all non marked boundary components of $\Sigma$. Each element of $\cH$ is a central element of $\cS\SM$, and the polynomial algebra $\BC[\cH]$ embeds into $\cS\SM$ because the monomials in elements of $\cH$ are part of the basis $B$ given by Theorem \ref{thm.basis}.

For each marked boundary component $\beta$ of $\Sigma$ consider the product
$\prod_{e\subset \beta} e$ of all boundary $\cP$-arcs on $\beta$. These stated arcs $q$-commute, so the product is well defined up to a scalar. Let $\mathcal{G}$ be the set of such elements.

The following theorem says that the center of the skein algebra, when $\xi$ is not a root of 1, is the obvious one.

\begin{theorem} 
Suppose $\xi$ is not a root of unity.

(a) \cite{LY} The center of $\stateS_\xi(\Sigma,\marked)$ is $\BC[\cH]$.

(b) \cite{LP} The center of $\stateS^+_\xi(\Sigma,\marked)$ is $\BC[\cH,\mathcal{G}]$.
\end{theorem}
When there is no marked point, part (a) was already proved in \cite{PS2}.

\subsection{Roots of 1: dimension} \label{sec.dim}
We already saw that when $\xi$ is not a root of 1, the center of $\cS\SM$ is small. At a root of unity, the center of $\cS\SM$ is considerably bigger. From the orderly finite generation of Theorem \ref{thm.basis}(b) and Corollary \ref{r.center1} one can 
easily show that $\cS_\xi\SM$ is a finitely-generated module over its center. Hence Theorem \ref{thm.Azu} applies to $\cS_\xi\SM$. It is important to calculate the PI dimension of $\cS_\xi\SM$, its center, and the Azumaya locus.



Recall that $r= 3m - 3\chi$ is the Gelfand-Kirillov dimension of $\cS_\xi\SM$, where $m=|\cP|$ is the number of marked points and $\chi$ is the Euler characteristic of $\Sigma$.
\begin{theorem}\label{thm-dim-z}
Let $\SM$ be a connected marked surface of genus $g$ and $\xi$ be a root of 1 with $N=\ord(\xi^4)$. The dimension of $\stateS(\Sigma,\marked)$ over its center is
\[D_\xi\SM=\begin{cases}
N^E,&\text{if }\ord(\xi) \text{ is odd},   \\
2^{2\lfloor\frac{m'}{2}\rfloor}N^E,&\text{if } \ord(\xi) \equiv 2 \pmod 4,\\
2^{2g+2m'}N^E,&\text{if }\ord(\xi) \equiv 0 \pmod 4.
\end{cases}\]
Here $m'= \max(0,m-1)$ and $E=r- b_\ev=3|\cP| - 3\chi(\Sigma)-b_\ev$, where $b_\ev$ is the number of boundary components of $\Sigma$ having an even number (including 0) of marked points.
\end{theorem}

One can show that $E$ is the rank of the matrix $Q$ in Theorem \ref{thm.embedding} and \ref{thm.embedding2} (or rank of $P$ in Theorem \ref{thm.embedding} whenever $P$ is defined), which means $b_\ev$ is the nullity of $Q$ (or $P$). By PI algebra theory, $E$ is even since $D_\xi$ is a perfect square, though this is not immediately obvious from $E=r- b_\ev$.

When $\SM$ has no marked points, Theorem \ref{thm-dim-z} was first proved in \cite{FKL2}. 

\subsection{Roots of 1: center}
Let $\SM$  be a marked surface and $\xi$ be a root of unity of order $n$, with  $N=\order(\xi^4)$ and $\epsilon=\xi^{N^2}$. The following are some special central elements.

\begin{lemma}\label{lemma-central-boundary}
Let $\xi$ be a root of unity of order $n$, with  $N=\order(\xi^4)$. Suppose $e_1,\dots,e_s$ are the boundary $\cP$-arcs on a boundary component $\beta$ with $s$ marked points, ordered consecutively with an arbitrary starting arc.
\begin{enumerate}
\item If $s$ is odd, the element $e_1^Ne_2^N\dots e_s^N$ is central.
\item If $s$ is even, the element $\hat{e}_1^k\hat{e}_2^{n-k}\dots \hat{e}_{s-1}^{k}\hat{e}_s^{n-k}$ is central for $0\le k\le n$.
\end{enumerate}
\end{lemma}

A tangle diagram on $\Sigma$ with endpoints at $\cP$ on $\SM$ is \emph{matching} if the parity of the number of end points at every marked point is the same. The {\em matching subalgebra} $\matchedS(\Sigma,\marked)$ is the span of matching elements.

A stated tangle diagram $\alpha$ is \emph{even} if $\alpha$ represents $0$ in $H_1(\bar{\Sigma},\marked;\mathbb{Z}/2)$, and at each marked point, the difference between the number of $+$ and $-$ states of $\alpha$ is divisible by $4$. Here $\bar{\Sigma}$ is the surface obtained by capping off all unmarked boundary components. The {\em even subalgebra} $\evenS(\Sigma,\marked)$ is the span of even elements.

From Theorem \ref{r.transparent} it is not difficult to show that the subalgebra
\begin{equation}
Z_0:=\begin{cases}
\Phi_\xi(\stateS_\epsilon(\Sigma,\marked))[\cH],& \text{if $\ord(\xi)$ is odd},\\
\Phi_\xi(\matchedS_\epsilon(\Sigma,\marked))[\cH],& \text{if }\ord(\xi) \equiv 2 \pmod 4,\\
\Phi_\xi(\evenS_\epsilon(\Sigma,\marked))[\cH],& \text{if }\ord(\xi) \equiv 0 \pmod 4
\end{cases}
\end{equation}
is in the center of $\cS_\xi\SM$.

\begin{theorem}\label{thm-center-root-1}  The center $Z(\stateS_\xi(\Sigma,\marked))$ of $\stateS_\xi(\Sigma,\marked)$  is $Z_0[C_\xi^{-1}]$, where $C_\xi$ is the set of central elements given by Lemma~\ref{lemma-central-boundary}. The affine variety $\MaxSpec(Z(\stateS_\xi(\Sigma,\marked)))$ is a normal variety of  dimension $r\SM= 3 |\cP|- 3 \chi(\Sigma)$.
\end{theorem}

Elements in $C_\xi$ are not invertible in $\cS\SM$. Instead, what the theorem means is that the center is spanned by elements $x$ such that $xc $ is in $Z_0$ for some product $c$ of elements in $C_\xi$.

When $\SM$ has no marked point, there is no $C_\xi$ and Theorem \ref{thm-center-root-1} simply says that the center $Z(\stateS_\xi(\Sigma,\marked))$ is $Z_0$. This case was proved in \cite{FKL}.

\def\Tr{\mathrm{Tr}}
\subsection{Cayley-Hamilton algebra, maximal order, Poisson order} Let $A$ be a $\BC$-algebra with assumptions as in Subsection \ref{sec.Azu}. By left multiplication, every $a\in A$ acts $\tZ$-linearly on $\tA$, a finite-dimensional $\tZ$-vector space. Denote by $\Tr(a)\in \tZ$ the trace of this action. If $\Tr(a)\in Z$ for all $a\in A$, then $A$ belongs to the class of {\em Cayley-Hamilton algebras} in the sense of \cite{DP}. For this class of algebras, there are finer results about the Azumaya locus than Theorem \ref{thm.Azu}. For example \cite{BY} shows that the Azumaya locus is the complement of the zero set of the determinant ideal, defined there.

The class of {\em maximal orders} \cite{MR,DP} is even finer than the class of Cayley-Hamilton algebras. Examples of maximal orders are quantum groups at roots of 1, see e.g. \cite{DP,BG}.
\begin{theorem} \cite{LY} \label{thm.tr} 
Let $\SM$ be a marked surface and $\xi$ be a root of 1.

(a) For every  $a\in \cS_\xi\SM$ one has $\Tr(a) \in Z(\cS\SM)$. In particular, $\cS_\xi\SM$ is Cayley-Hamilton and the Azumaya locus is the complement of the zero locus of the determinant ideal.

(b) If $\SM$ is not a closed surface $\Sigma_g$, then $\cS\SM$ is a maximal order.
\end{theorem}

For the case when $\SM$ has no marked points, part (a) was proved in \cite{FKL2}, with partial results in \cite{AF}. Also in this case, part (b) was proved by the first author and J. Paprocki, see \cite{Paprocki}. Let us sketch a proof of (b) for general marked surfaces. It is known that \cite{MR} if the associated graded algebra of $A$ with respect to an $\BN$-filtration is a maximal order, then  $A$ is a maximal order. Theorem \ref{thm.embedding}(b) says an associated graded algebra of $\cS\SM$ with respect to an $\BN$-filtration is a monomial algebra $A(Q;\Lambda)$ for a certain submonoid $\Lambda\subset \BZ^r$. Explicit calculation shows that this submonoid is primitive in the sense that if $k\lambda\in \Lambda$ where $\lambda\in \BZ^r$ and $k$ is a positive integer, then $\lambda\in \Lambda$. It then follows from  \cite{Paprocki} that $A(Q;\Lambda)$ is a maximal order at roots of 1. 

Besides, in all cases, $\cS_\xi\SM$ is a Poisson order over its center in the sense of \cite{BG}, which allows us to use the Poisson geometry of its center to study the Azumaya locus.

\subsection{Azumaya locus: Punctured torus} For the punctured torus we have a full description of the Azumaya locus.

\begin{theorem}\cite{LY}
Let $\Sigma_{1,1}$ be the compact oriented surface of genus 1 and having one boundary component and no marked points, and $\xi$ be a root of 1. The Azumaya locus of $\cS_\xi(\Sigma_{1,1})$ is the set of smooth (or regular) points of the affine variety $\MaxSpec(Z(\cS_\xi(\Sigma_{1,1})))$, which  is a finite branched covering of $\BC^3$ and has a finite number of non-smooth points.
\end{theorem}

In \cite{LY} we give an explicit description of the branch covering and the singular set.

\subsection{Azumaya locus: Closed surfaces}
Recall that $\Sigma_g$ is the closed surface of genus $g$.

\begin{theorem} Let $\xi$ be a root of 1 and $g\ge 2$.

(a) The Azumaya locus of $\cS(\Sigma_g)$ contains the smooth locus of
$\MaxSpec(Z(\cS_\xi(\Sigma_g)))$. 

(b) The variety $\MaxSpec(Z(\cS_\xi(\Sigma_g)))$ is isomorphic to the $SL_2(\BC)$-character variety $\cX(\Sigma_g)$ if $\ord(\xi) \neq 0 \pmod 4$, and is isomorphic to the connected component of the $PSL_2(\BC)$-character variety which contains the character of the trivial representation, if $\ord(\xi) =0 \pmod 4$.
\end{theorem}

When $\ord(\xi)$ is odd, part (a) was proved in \cite{Jordan}. The extension to all roots of 1 will be considered in \cite{LY}. Part (b) was proved in \cite{FKL}. It is known that the smooth part of $\cX(\Sigma_g)$ consists exactly of the characters of irreducible $SL_2(\BC)$-representations of the $\pi_1(\Sigma_g)$.

We conjecture that the Azumaya locus is equal to the smooth locus, for $\Sigma_g$.


\end{document}